%% file: unification.tex
\documentclass[11pt]{article}
\usepackage{graphicx} 
\usepackage{amsmath,amssymb,latexsym,proof}

\input{macros}

\newcommand{\nmodels}{\nvDash}
\newcommand{\GLB}{\mathrm{GLB}}
\renewcommand{\GLP}{\mathrm{GLP}}
\renewcommand{\GL}{\mathrm{GL}}
\renewcommand{\J}{\mathrm{J}}
\newcommand{\Fm}{\mathrm{Fm}}
\newcommand{\GLPz}{\mathrm{GLP.3}}
\renewcommand{\PA}{\mathrm{PA}}

\title{On the unification problem for GLP}
\author{L.D. Beklemishev}

\begin{document}

\maketitle

\begin{abstract}
    We show that the polymodal provability logic $\GLP$, in a language with at least two modalities and one variable, has nullary unification type. More specifically, we show that the formula $[1]p$ does not have maximal unifiers, and exhibit an infinite complete set of unifiers for it. Further, we discuss the algorithmic problem of whether a given formula is unifiable in $\GLP$ and remark that this problem has a positive solution. Finally, we state the arithmetical analogues of the unification and admissibility problems for $\GLP$ and formulate a number of open questions.
\end{abstract}

\section{Introduction}

We work in the language $\cL(n,m)$ of propositional polymodal logic with modalities traditionally denoted $[0],[1],\dots , [n-1]$ and variables $\{p_0,\dots,p_{m-1}\}$. Let $\Fm(n,m)$ denote the set of all formulas of $\cL(n,m)$. A \emph{substitution} is a map $\gs: \Fm(n,m) \to \Fm(n,m)$ preserving all logical connectives and modalities. A \emph{variable-free substitution} is a map $\gs: \Fm(n,m) \to \Fm(n,0)$.

By a \emph{logic} here we mean a set of formulas closed under the axioms and rules of classical propositional logic and the substitution rule: from $\phi$ infer $\gs(\phi)$. \emph{Normal logics} are those containing the axioms of basic modal logic K for each modality and closed under modus ponens and necessitation rules: from $\phi$ infer $[i]\phi$. If $\Gamma$ is a finite set of formulas, we write $\Gamma\vdash_L \phi$ if there is a proof of $\phi$ from the axioms of $L$ and $\Gamma$ using modus ponens and necessitation rules.

\bd
A \emph{unifier} for $\phi$ in a logic $L$ is a substitution $\gs$ such that $L\vdash \gs(\phi)$.

A unifier $\gs$ is \emph{more general} than $\tau$ in $L$ (denoted $\tau\leq \gs$) if there is a substitution $\theta$ such that $\tau=\theta\gs$ in $L$, that is, for all variables $p$,
$$L\vdash \tau(p)\eqv \theta(\gs(p)).$$
\ed

Let $U_L(\phi)$ denote the set of all unifiers of a formula $\phi$ in $L$. $\phi$ is called \emph{unifiable} if $U(\phi)$ is not empty. A unifier $\gs\in U(\phi)$ is called \emph{maximal} for $\phi$ if there is no $\tau\in U(\phi)$ such that $\gs < \tau$. Here, as usual, $\gs<\tau$ means $\gs\leq\tau$ and not $\tau\leq \gs$.

\bd A subset
$\Sigma\subseteq U_L(\phi)$ is a \emph{complete set of unifiers for $\phi$} if  $$\al{\tau\in U_L(\phi)}\ex{\sigma\in\Sigma} \tau\leq\sigma.$$
\ed

\bd
A logic $L$ has \emph{finitary/infinitary  unification type} if, for each unifiable formula $\phi$, there is a finite/infinite complete set of maximal unifiers for $\phi$.
$L$ has \emph{nullary unification type} if it is neither finitary nor infinitary. In this case, there is a unifiable formula $\phi$ without a complete set of maximal unifiers.
\ed

Silvio Ghilardi in the 1990s (see~\cite{Ghi99,Ghi00}) pioneered the study of unification in propositional logics; his methods were based on the notion of projective formula and projective substitution. He showed that many standard logics, including the intuitionistic logic IPC and the modal logics K4, S4, GL, Grz, have finitary unification type, and complete sets of maximal unifiers can be constructed effectively. As a corollary, this yielded another solution of the problem of \emph{admissibility of inference rules} for these logics previously obtained by Vladimir Rybakov by different methods (see~\cite{Ryb84,Ryb89,Ryb97}). Explicit axiomatizations of admissible rules based on (extensions of) Ghilardi's methods for various logics have been obtained by R. Iemhoff, E. Je\v{r}\'{a}bek, G. Metcalfe, V. Rybakov, and many others (see \cite{Iem01,Jera05,IemMet09,Bezh16}).

Ghilardi's methods essentially rely on \emph{transitivity} and \emph{finite model property} (FMP) of the considered logics. For the basic (intransitive) modal logic K the decidability of the unification problem and of the admissibility of rules problem are still open. E. Je\v{r}\'{a}bek~\cite{Jera15} has shown that K has nullary unification type.

For polymodal logics the situation is even more complicated. For example, Wolter and Zakharyaschev~\cite{WoZa08} show that for the bimodal logics extending K and K4 by the universal modality the unification problem and the admissibility of rules problem are undecidable. In contrast, the derivability problem for these logics is well-known to be decidable.

In this paper, we deal with the polymodal provability logic $\GLP$ that is transitive, but does not have FMP (nor is Kripke complete).
The logic GLP introduced by G. Japaridze~\cite{Dzh86} plays an important role in proof-theoretic applications of provability logic. It has also been the subject of rather detailed modal-logical studies over the last years. In particular, we know that $\GLP$ satisfies Craig's interpolation property~\cite{Ign93}, as well as uniform interpolation~\cite{Sham11}. Its variable-free fragment has effective normal forms~\cite{Ign93,Bek04}, and the elementary theory of the Lindenbaum algebra of the variable-free fragment of $\GLP$ (in the language with $n$ modalities) is decidable~\cite{Pakh16}.

The unification and the admissibility of rules problems for $\GLP$ have been stated as open questions in~\cite{BekJoo12}. The first steps in the study of the unification problem for $\GLP$ were undertaken by N. Lukashov in his Bachelor's Thesis~\cite{LukBT}. He considered the subsystem $\J_2$ of $\GLP_2$ and proved that $\J_2$ has finitary unification type. It is known that the derivability problem for $\GLP_2$ is reducible to $\J_2$, however, as we show here, the reduction does not preserve the unification type.

We show that the fragment $\GLP_n$ of $\GLP$ in the language with $n$ modalities has nullary unification type.
More specifically, the formula $[1]p$ does not have maximal unifiers, and we exhibit an infinite complete set of unifiers for it. Furthermore, we discuss the algorithmic problem of whether a given formula is unifiable in $\GLP_n$ and  remark that this problem has a positive solution, which can be inferred from the results of F.~Pakhomov and J.~Aguilera. In the last section of the paper we point out a few interesting connections of the unification problem for $\GLP$ with some recent results in the proof-theory of arithmetic.

\section{Japaridze's provability logic GLP}
\paragraph{Logics $\GLP$ and $\J$.}
The polymodal provability logic $\GLP$ is formulated in the language
of the propositional calculus enriched by new unary connectives
$[0]$, $[1]$, $[2]$, \dots , called \emph{modalities}. The dual
connectives $\la n\ra$, for all $n\in\nat$, are treated as
abbreviations: $\la n\ra\phi$ means $\neg [n]\neg\phi$.

The logic $\GLP$ is given by the following axiom schemata and
inference rules.
\begin{description}
\item[Axiom schemata:]
\begin{enumr}
\item Tautologies of classical propositional calculus; \item $[n](\phi\imp\psi)\imp
([n]\phi\imp [n]\psi)$; \item $[n]([n]\phi\imp \phi)\imp [n]\phi$
(L\"ob's axiom);
\item $\la m\ra\phi
\imp [n]\la m\ra\phi$, for $m<n$;
\item $[m]\phi\imp [n]\phi$, for
$m\leq n$ (monotonicity axiom);
\end{enumr}
\item[Inference rules:] modus ponens, $\vdash\phi\ \Imp\ \vdash [n]\phi$.
\end{description}

Thus, for each modality $[n]$ $\GLP$ contains the axioms
(i)--(iii) and the inference rules of G\"odel--L\"ob provability
logic $\GL$, whereas the schemata (iv) and (v) tie together
different modalities. A well-known consequence of axioms
(i)--(iii) is the \emph{transitivity schema}
$$[n]\phi\to [n][n]\phi.$$

Logic $\J$ is obtained from $\GLP$ by replacing the monotonicity
axiom (v) by the following two axiom schemata easily derivable in
$\GLP$ using the transitivity schema:

\begin{enumr} \setcounter{enumi}{5}
\item $[m]\phi\to [n][m]\phi$, for $m<n$;
\item $[m]\phi\to [m][n]\phi$, for $m<n$.
\end{enumr}

When the language is restricted to $n$ modalities $[0],\dots,[n-1]$, the systems axiomatized by the rules of
$\GLP$ and $\J$ are denoted $\GLP_n$ and $\J_n$, respectively. Thus, $\GLP_1$ is the G\"odel-L\"ob logic $\GL$.
It is known that $\GLP$ is conservative over $\GLP_n$, $\GLP_n$ is conservative over $\GLP_m$, for $m<n$, and similar relations hold for $\J$ and $\J_n$.

Unlike $\GLP$, the logic $\J$ is complete w.r.t.\ a natural class of
finite Kripke frames described below.

\paragraph{Kripke frames.}
A \emph{Kripke frame} for the language of polymodal logic is a
structure $(W;R_0,R_1,\ldots)$, where $W$ is a nonempty set, and
$R_k$ are binary relations on $W$. The elements of $W$ are usually
called the \emph{worlds}. A frame is called finite if $W$ is finite
and $R_k=\emptyset$, for all but finitely many $k\geq 0$.

A \emph{valuation} $v$ on a frame $W$ maps every propositional
variable $p$ to a subset $v(p)\subseteq W$ called the \emph{truth
set of $p$}. A \emph{Kripke model} is a Kripke frame together with a
valuation on it.

Let $\cW=(W,v)$ be a Kripke model. By induction on the build-up of
$\phi$ we define a relation \emph{$\phi$ is true at $x$ in $\cW$} (denoted $\cW,x\models\phi$).
\ben
\item $\cW,x\models p \iff x\in v(p)$, if $p$ is a variable;
\item $\cW,x\models (\phi\land\psi) \iff (\cW,x\models \phi\text{ and }\cW,x\models \psi)$, \\[1ex]
$\cW,x\models \neg\phi \iff \cW,x\nmodels \phi$, \\[1ex]  and similarly for the other boolean connectives;
\item $\cW,x\models [n]\phi \iff \al{y\in W}(xR_n y \Imp \cW,y\models\phi)$.
\een

Write $\cW\models\phi$ if $\al{x\in W} \cW,x\models\phi$.

For the axioms of $\J$ to be valid in a given Kripke frame for any
valuation of variables, we impose some restrictions on the relations
$R_k$. A transitive binary relation $R$ on a set $W$ is called
\emph{noetherian}, if there is no infinite chain of elements of $W$
of the form $a_0Ra_1Ra_2\dots$ Note that if $W$ if finite the
condition of $R$ being noetherian on $W$ is equivalent to its
irreflexivity.

A frame $(W;R_0,R_1,\ldots)$ is called a \emph{$\J$-frame} if

\ben
\item $R_k$ is transitive and noetherian, for all $k\geq 0$;
\item $\forall x,y\:( xR_ny \Imp \forall z\:(xR_mz \Leftrightarrow
yR_mz))$, for $m<n$;
\item $\al{x,y,z}(xR_m y \text{ and } yR_n z \Imp xR_m z)$, for $m<n$.
\een

A \emph{$\J$-model} is any Kripke model based on a $\J$-frame. By
induction on the derivation length one can easily prove the
following lemma.

\bl For any formula $\phi$, if $\J\vdash \phi$ then
$\cW\models\phi$, for any $\J$-model $\cW$. \el

The converse statement, that is, the completeness theorem for $\J$
with respect to the class of all (finite) $\J$-models, is proved in
ref.~\cite{Bek10}. Let us call a \emph{root} of a $\J$-frame $(W;R_0,R_1,\dots)$ an
element $r\in W$ for which $$\al{x\in W} \ex{k\geq 0} r R_k x.$$

\bpr \label{jcomp1} For any formula $\phi$, if $\J\nvdash\phi$ then
there is a finite $\J$-model $\cW$ with a root $r$ such that
$\cW,r\nmodels\phi$. \epr

Define  $$M^+(\phi)=M(\phi)\land \bigwedge_{i\leq n}[i]M(\phi),$$ where $n$ is the maximum modality number occurring in $\phi$, and $M(\phi)$ is the conjunction of all formulas of the form $[i]\psi \to[j]\psi$ with $[i]\psi\in \mathit{Sub}(\phi)$ and $i<j\leq n$.
A relationship between $\J$ and $\GLP$ is stated in the following proposition (see~\cite[Theorem 4]{Bek10}).
\bpr \label{redj} $\GLP\vdash \phi$ iff $\J\vdash M^+(\phi)\to \phi$.
\epr

In particular, this means that the derivability problem for $\GLP$ is reducible to the one for $\J$.

\section{Unifiers of $[1]p$} \label{boxone}

Firstly, we notice that the formula $[1]p$ has a unifier $p\mapsto ([0]p\to p)$ in $\GLP$. Curiously, this fact corresponds (under the standard provability logical interpretation) to a well-known theorem  by S. Feferman~\cite{Fef60} stating that the theory of all true $\Pi_1$-sentences in Peano arithmetic proves its local reflection schema. The corresponding  modal logical lemma is just an abstraction of his proof.

\bl \label{one} $\GLP\vdash [1]([0]q\to q)$. \el

\bp\ On the one hand, $\GLP$ proves that
$[0]q$ implies $[1]q$ and hence $[1]([0]q\to q)$. On the other hand, $\neg [0]q$ implies $[1](\neg [0]q)$ and hence $[1]([0]q\to q)$. In each case we obtain the required.
\ep

\bcor \label{three}
For each $k$ the substitution $p \mapsto \langle 0\rangle^k \top$ is a unifier for $[1]p$ in $\GLP$.
\ecor

\bp\
This follows by induction on $k$ using that $\GLP\{[0]q\to q\}\vdash [0]^{m+1} \bot \to [0]^m\bot$, for all $m$.
\ep

The following lemma characterizes all unifiers of the formula $[1]p$ in $\GLP$, and thereby in each of the logics $\GLP_n$, regardless of the number of variables. We notice that each unifier corresponds to a formula $\phi$ such that $\GLB\vdash [1]\phi$.

\bd
Let $\GLP\{[0]q\to q\}$ denote the (semi-normal) logic obtained by the closure under \emph{modus ponens} and the substitution rules of the set of all theorems of $\GLP$ and the formula $[0]q\to q$.
\ed

\bl \label{two} $\GLP\vdash [1]\phi$ iff $\GLP\{[0]q\to q\}\vdash \phi$.
\el

\bp\ The (if) part is clear, since the set of all $\phi$ such that $\GLP\vdash[1]\phi$ contains all theorems of $\GLP$, the formula $[0]q\to q$ (by Lemma \ref{one}) and is closed under \emph{modus ponens} and the substitution rules.

To prove the (only if) part,
assume that  $\GLP\{[0]q\to q\}\nvdash \phi$. Let $R(\phi)$ denote $$\bigwedge_{\psi \in Sub(\phi)}([0]\psi \to \psi),$$ then $\GLP \nvdash R(\phi) \to \phi$.
Since $\GLP\vdash M^+([1]\phi)$, we have $$\J \nvdash M^+([1]\phi)\land R(\phi) \to \phi.$$
Notice that we put $M^+([1]\phi)$ rather than $M^+(\phi)=M^+(R(\phi)\to \phi)$ here.
There is a $\J$-model $\mathcal{K}$ with a root $r$ such that $\mathcal{K}, r \models M^+([1]\phi), R(\phi)$ and $\mathcal{K}, r \nmodels \phi$. In particular, at every node of $\mathcal{K}$ the formula $M([1]\phi)$ is true. We want to transform $\mathcal{K}$ to a $\J$-model falsifying  $M^+([1]\phi)\to [1]\phi$.

Let $A$ denote the subset $\{x\in \cK: r R_1 x\}\cup \{r\}$.
Let a model $\mathcal{K}'$ be obtained from $\mathcal{K}$ by adding a new root $r'$ to the restriction of $\mathcal{K}$ to $A$. In other words, let $r' R_1 x$, for all $x\in A$, and $r' R_0 y$, for all $y$ such that $r R_0 y$. For $i>1$, $r' R_i z$ holds for no $z$.  The evaluation of variables at $r'$ is arbitrary. (We consider $\cK'$ as an extension of $\cK$ and do not rename the relations.) Obviously, the truth value of all formulas at the nodes of $\cK$ is preserved in $\cK'$.
Since $r' R_1 r$ and $\cK',r\nmodels \phi$, we have $\cK',r'\nmodels [1]\phi$.

We claim that $\mathcal{K}', r' \models M^+([1]\phi)$. It is sufficient to check that $\mathcal{K}', r' \models M([1]\phi)$, since at all the nodes from $\mathcal{K}$ the formula $M([1]\phi)$ remains true in $\cK'$.

Let $[i]\psi$ be a subformula of $[1]\phi$. Clearly, the formulas $[j]\psi$, for all $j>1$, are true at $r'$, hence $[i]\psi\to [j]\psi$ is also true. So, we only need to consider $i=0$, that is, to show that $\mathcal{K'},r'\models [0]\psi\to [1]\psi$, for all subformulas $[0]\psi$ of $\phi$.

Assume $[0]\psi \in \mathit{Sub}(\phi)$ and $\mathcal{K}', r' \models [0]\psi$. For any $x$ such that $r' R_1 x$ we need to show that $\cK',x\models \psi$; this would entail the required $\cK',r'\models [1]\psi$.

First, consider $x=r$. Since the same points are $R_0$-accessible from $r$ as from $r'$, it follows that $\mathcal{K}', r \models [0]\psi$ and $\mathcal{K}', r \models \psi$ (because $\mathcal{K}, r \models R(\phi)$).

Suppose $x\in A\setminus \{r\}$. Since $\mathcal{K}', r \models [0]\psi$ and $M(\phi)$ holds at $r$ we obtain  $\mathcal{K}', r \models [1]\psi$. But $r R_1 x$, so we obtain $\mathcal{K}', x \models \psi$. This concludes the proof that $\mathcal{K}', r' \models M^+(\phi)$ and obviously $\mathcal{K}', r' \nmodels [1]\phi$.
\ep

Define a sequence of formulas $Q_1(p):=p$ and $Q_{i+1}(p):=(p\lor [0] Q_i(p))$ and let $Q^n(p):=\bigwedge_{i=1}^n([0]Q_i\to Q_i).$
We will rely on the following lemma which is a slight stengthening of~\cite[Lemma 2.27]{Bek05en}.

\bl \label{main} $\bigwedge_{i=1}^n ([0]p_i\to p_i) \to p \vdash_{\GL} Q^n(p)\eqv  p$.
\el
Informally, this lemma is saying that if $p$ follows from \emph{any} $n$ instances of reflection, then $p$ is equivalent to a conjunction of very specific $n$ instances of reflection (dependent only on $p$). For reader's convenience we give a short proof here.

\bp\
First, we remark that $\GL\vdash p\to Q_i(p)$, for all $i$, hence $\GL\vdash p \to Q^n(p)$. For the other part of the equivalence
we establish two deductions: $$\bigwedge_{i=1}^n ([0]p_i\to p_i) \to p \vdash_\GL Q_{n+1}(p) \vdash_\GL Q^n(p)\to  p.$$

Part 1. We use a Kripke model argument. Consider any Kripke model $\cW=(W,R_0,v)$ and an $x_0\in W$ such that $\cW,x_0\nmodels Q_{n+1}(p)$. Then there is an increasing chain of nodes $x_0 R_0x_1 R_0\dots R_0 x_n$ such that $\cW,x_k\nmodels p$, for all $k$. Each $[0]p_i\to p_i$ is false in at most one node of the chain. By the pigeon-hole principle, there is a $k$ such that LHS (denote this formula $A$) is true at $x_k$. Then, $\cW,x_k\nmodels A\to p$ and hence $\cW,x_0\nmodels [0]^+(A\to p)$. Since the model $\cW$ was arbitrary, this shows $\GL\vdash [0]^+(A\to p) \to Q_{n+1}(p)$.

Part 2. By a simple induction on $i\leq n+1$ we prove that
$$Q^n,Q_{n+1}\vdash p\lor Q_{n+1-i}.$$
For $i=0$ the claim is trivial. Otherwise we use the induction hypothesis and infer:
$$Q^n, p\lor Q_{n+1-i}\vdash p\lor [0]Q_{n+1-i-1}\vdash p\lor Q_{n-i}.$$
Putting $i=n$ entails the claim.
\ep

Somewhat abusing the language, we will call a formula $\phi$ such that $\GLP\vdash[1]\phi$ a \emph{unifier of $[1]p$}. (We identify $\phi$ with the substitution $p\mapsto \phi$.) From Lemma~\ref{main} we obtain the following corollary.

\bt A substitution $\sigma$ is a unifier of $[1]p$ in $\GLP$ iff $\sigma\leq Q^k$, for some $k$.
\et
\bp\ By Lemma~\ref{two} all $Q^k$ are unifiers of $[1]p$. Conversely, if a formula $\psi=\sigma(p)$ is a unifier of $[1]p$, then by Lemma~\ref{two} $\psi$ follows from some $k$ instances of reflection in $\GLP$. By Lemma \ref{main} (with $\psi$ substituted for $p$) we obtain $\GLP\vdash \psi \eqv Q^k(\psi)$, which means $\GLP\vdash \gs(p)\eqv \gs(Q^k(p))$, that is, $\sigma\leq Q^k$. \ep

This proof shows a bit more:

\bcor $Q^k$ is the most general unifier for $[1]p$ among those which follow from $\leq k$ instances of reflection. \ecor

In particular, we have $\la 0\ra^{n}\top\leq Q^n$. It is well-known by the results of Boolos and Artemov that $\la 0\ra^{n}\top$ does not follow from less than $n$ instances of reflection in $\GL$. The same result holds for $\GLP$, as the following lemma shows.

\bl \label{abl} For any formulas $\psi_1,\dots, \psi_n$, $$\bigwedge_{i=1}^n ([0]\psi_i\to \psi_i) \nvdash_\GLP \la 0\ra^{n+1}\top.$$ \el

\bp From Lemma~\ref{main}, substituting $\bot$ for $p$, we obtain the following corollary (see~\cite{Bek05en}):
\begin{equation} \GL\vdash \la0\ra^{n+1}\top \to \la 0\ra\bigwedge_{i=1}^n ([0]p_i\to p_i).\label{ab} \end{equation}
Now assume, for a contradiction, that $$\bigwedge_{i=1}^n ([0]\psi_i\to \psi_i) \vdash_\GLP \la 0\ra^{n+1}\top.$$
Then also $$\GLP\vdash \la 0\ra\bigwedge_{i=1}^n ([0]\psi_i\to \psi_i) \to \la 0\ra^{n+2}\top.$$
Therefore, by \refeq{ab} we obtain $\GLP\vdash \la0\ra^{n+1}\top \to \la 0\ra^{n+2}\top$ whence by L\"ob's axiom $\GLP\vdash [0]^{n+1}\bot$, which is not the case.
\ep

\bcor For all $i,j>0$, $Q^i<Q^j$ whenever $i<j$. \ecor

\bp If $i<j$ then $Q^i\leq Q^j$, since $Q^i$ follows from $i<j$ instances of reflection. For the same reason, $\la 0\ra^{j}\top\leq Q^j$, yet by Lemma~\ref{abl} $\la 0\ra^{j}\top\nleq Q^i$. Otherwise, $\GLP\vdash \la 0\ra^{j}\top\eqv Q^i(\phi)$, for some $\phi$, and
 $\la 0\ra^{j}\top$ would follow from $i$ instances of reflection $Q^i(\phi)$ contradicting the lemma. \ep

\bcor For any $n>1$, the formula $[1]p$ has no maximal unifiers in $\GLP_n$. The unification type of $\GLP_n$ is nullary. \ecor

The same corollary also holds for the full language of $\GLP$ (with $\gw$-many modalities): a unifier of $[1]p$ in $\GLP$ is a unifier in $\GLP_n$, for some $n$, and less general than some $Q^k$. On the other hand, since $\GLP_1$ is just a notational variant of $\GL$, its unification type is finitary.

\section{Unifiers of $[1]p$ and the reduction property}
A sightly simpler characterization of the set of unifiers of $[1]p$ in $\GLP$ can be obtained via the so-called \emph{reduction property}.
The reduction property for free $\GLP$-algebras is established in~\cite{Bek17}, its arithmetical analogues play a significant role in the applications of $\GLP$. We state here without proof a somewhat less general proposition that is a particular case of Lemma 6 of~\cite{Bek17}.

Define a sequence of formulas $Q^*_n(p)$ as follows: $Q^*_1(p):=p$, $Q^*_{n+1}(p):=p\land \la 0\ra Q^*_n(p)$. Clearly, for each $n$, $Q^*_n(p)$ is equivalent to $\neg Q_n(\neg p)$ in $\GLP$.

\bpr \label{red}
Suppose $\GLP\vdash \la 1\ra\phi \to \la 0 \ra\psi$. Then there is a $k$ such that $\GLP\vdash Q^*_k(\phi)\to \la 0\ra\psi$. 
\epr

From this proposition we infer

\bl $\GLP\vdash [1]\phi$ iff $\GLP\vdash Q_k(\phi)$, for some $k$. \el

\bp Assume $\GLP\vdash [1]\phi$, then $\GLP\vdash \la 1\ra\neg\phi \to \bot$. Since $\bot$ is equivalent to $\la 0\ra\bot$, Proposition~\ref{red} is applicable and we obtain $\GLP\vdash Q^*_k(\neg\phi)\to \bot$, that is, $\GLP\vdash Q_k(\phi)$, for some $k$.

In the other direction, assume $\GLP\vdash Q_k(\phi)$, for some $k$. We prove $\GLP\vdash [1]Q_i(\phi)$, for all $1\leq i\leq k$, by downwards induction on $i$. For $i=k$ the claim is obvious, assume $i<k$ and the claim holds for $i+1$. Recall that $Q_{i+1}(\phi)=\phi\lor [0]Q_i(\phi)$. Assume $\GLP\vdash [1](\phi\lor [0]Q_i(\phi))$. Then, since $\GLP\vdash \phi\to Q_i(\phi)$, we have $\GLP\vdash [1](Q_i(\phi)\lor [0]Q_i(\phi))$.
By Lemma~\ref{one} we have $\GLP\vdash [1]([0]Q_i(\phi)\to Q_i(\phi))$, therefore $\GLP\vdash [1]Q_i(\phi)$. Putting $i=1$ we obtain $\GLP\vdash \phi$.
\ep

A useful form of this lemma is
\bcor $\GLP\vdash [1]\phi$ iff $\GLP\vdash \phi \eqv (Q_k(\phi)\to \phi)$, for some $k$. \ecor

\bp Clearly, $\GLP\vdash Q_k(\phi)$ holds iff $\GLP\vdash \phi \eqv (Q_k(\phi)\to \phi)$. \ep

\bcor A substitution $\gs$ is a unifier of $[1]p$ in $\GLP$ iff, for some $k$, $\gs\leq (Q_k(p)\to p)$. \ecor
\bp (only if) Putting $\phi:=\gs(p)$ in the above corollary we obtain that $\gs$ is a unifier of $[1]p$ iff $\GLP\vdash \gs(p)\eqv \gs(Q_k(p)\to p)$, for some $k$. The latter implies $\gs\leq (Q_k(p)\to p)$.

(if) By induction on $k$ we show that all formulas $Q_k(p)\to p$ are unifiers of $[1]p$. The statement clearly holds for $k=1$. To show $[1](Q_{k+1}(p)\to p)$ it is sufficient to derive $[1]([0]Q_k(p)\to p)$. By Lemma~\ref{one} we have $[1]([0]Q_k(p)\to Q_k(p))$, and by the induction hypothesis $[1](Q_k(p)\to p)$. Using the syllogism rule under $[1]$ we conclude $[1]([0]Q_k(p)\to p)$, as required. \ep

To round off the discussion, we mention without proof that the formulas $Q^n(p)$ are, in fact, equivalent to $Q_{n+1}(p)\to p$ in $\GLP$. Therefore, the two characterizations of the set of unifiers of $[1]p$ are essentially the same.

\section{The unification problem for $\GLP$}

The \emph{unification problem} for a logic $L$ is to determine, given a formula $\phi$, whether it is unifiable in $L$ or not. We show, in this section, that the unification problem for $\GLP$ and its fragments $\GLP_n$ is decidable.

First, observe that if $\phi$ is unifiable, then there is a unifier $\gs$ for $\phi$ in $L$ such that $\gs(p)$ is a variable-free formula, for each variable $p$ (for example, take any unifier $\theta$ and consider $\gs:=\tau\theta$ where $\tau$ substitutes $\top$ for each variable). Therefore, unifiability of a formula can be tested on the variable-free fragment of $L$.

The variable-free fragment of $\GLP$ and the associated Lindenbaum algebra are well-studied and play a prominent role in the proof-theoretic applications of $\GLP$, see~\cite{Ign93,BJV,Ica09}. Especially interesting for our purposes are two  recent papers.

Pakhomov~\cite{Pakh16} proves that the first order theory of the Lindenbaum algebra $\cA_n$ of the variable-free fragment of $\GLP_n$ is decidable. Unifiability of $\phi(p_1,\dots,p_m)$ (by a variable-free substitution) is equivalent to the truth of the existential statement $$\ex{x_1}\dots \ex{x_m} T_\phi(x_1,\dots, x_m)=1$$ in $\cA_n$, where $T_\phi$ denotes the term, corresponding to $\phi$. Thus, Pakhomov's theorem entails

\bt The unification problem for $\GLP_n$ is decidable. \et

The paper~\cite{Pakh16} is rather complicated and proves a much stronger result than what is really needed here. A simpler proof can be obtained on the basis of a recent paper by Aguilera and Pakhomov~\cite{AgPa24}. They show that the logic of variable-free substitutions of $\GLP$ can be axiomatized by a system $\GLPz$, a polymodal analogue of the logic $\mathrm{GL.3}$ of finite irreflexive linear frames. In other words, $\GLP\vdash\gs(\phi)$, for all variable-free substitutions $\gs$, if and only if $\GLPz\vdash\phi$. They also show the logic $\GLPz$ to be decidable. We do not state an explicit axiomatization of $\GLPz$ here, it is sufficient to use its characterization as the logic of variable-free substitutions in $\GLP$.

Our characterization of unifiable formulas in $\GLP_n$ is the following result.

\bt \label{unif} Let $\phi\in\cL(n,m)$. Then $\GLPz\nvdash \la n\ra\top \to \neg [0]\phi$ iff
there is a variable-free substitution $\gs:\cL(n,m)\to \cL(n,0)$ such that $\GLP_n\vdash \gs(\phi)$.
\et
Notice that this reduces the unifiability problem for $\GLP_n$ (in the language $\cL(n,m)$) to the derivability problem for $\GLPz$ in the language $\cL(n+1,m)$. Our proof is based on certain common techniques for the variable-free fragment of $\GLP$.

Recall that \emph{words} or \emph{worms} are variable-free formulas obtained from $\top$ by operations $\la n\ra$, for any $n<\gw$. It is well-known that the set of worms is closed under conjunction modulo provability in $\GLP$~\cite[Lemma 9]{Bek05a}.

\begin{lemma} \label{worms} If $\phi$ in $\cL(n,0)$ is variable-free, then $\la 0 \ra\phi$ is $\GLP_n$-equivalent to either $\bot$ or a worm of the form $\la 0\ra A$, for some $A$ in $\cL(n,0)$.
\end{lemma}
\bp\ We will infer this claim from some results in \cite{Bek05a} (other ways of proving it are surely possible). Bring $\phi$ to a disjunctive normal form. Since $\la 0\ra$ commutes with disjunction, $\la 0\ra\phi$ is equivalent to a disjunction of formulas of the form $\la 0\ra\psi$ where $\psi$ is a conjunction of worms and negated worms. Since conjunction of worms is equivalent to a worm, \cite[Lemma 10]{Bek05a} is applicable and yields that each of the formulas $\la 0\ra\psi$ is equivalent to $\bot$ or to a worm of the form $\la 0\ra A$. One can delete the disjuncts of the form $\bot$ as long as the disjunction is not empty. By \cite[Corollary 6]{Bek05a} worms of the form $\la 0 \ra A$ are linearly ordered by implication in $\GLP$. So, the disjunction is equivalent to the weakest one of them. \ep

\bl \label{vf} For any variable-free formula $\phi$ in $\cL(n,0)$, $$\GLP_n\nvdash \phi \iff \GLP\vdash \la n\ra\top\to\neg [0]\phi.$$
\el

\bp  It is well-known that $\GLP_n\nvdash\phi$ iff $\GLP_n\nvdash [0]\phi$ iff $\GLP_n\nvdash \neg[0]\phi\eqv\bot$. By Lemma~\ref{worms}, the latter holds if and only if $\GLP_n\vdash \neg[0]\phi \eqv A$, for some worm $A\in\cL(n,0)$. However, for any such $A$ we have $\GLP\vdash \la n\ra\top \to A$, which is readily seen by induction on the length of $A$. Thus, if $\GLP_n\nvdash\phi$ then $\GLP\vdash \la n\ra\top\to \neg[0]\phi$.

In the other direction, if $\GLP\vdash \la n\ra\top \to \neg[0]\phi$ then $\neg[0]\phi$ cannot be equivalent to $\bot$ (since $\GLP\nvdash [n]\bot$), therefore $\GLP_n\nvdash\phi$. \ep

Suppose $\phi$ is a variable-free formula of $\cL(n+1,0)$. Let $\ol\phi$ denote the result of replacing every subformula of $\phi$ of the form $[n]\psi$ by $\top$. Obviously, $\ol\phi\in\cL(n,0)$.

\bl \label{bn} Suppose $\GLP\nvdash \la n\ra\top\to \neg [0]\phi$, then $\GLP\nvdash \la n\ra\top\to \neg [0]\ol{\phi}$.
\el

\bp If $\GLP\nvdash \la n\ra\top\to \neg [0]\phi$, then there is a node $\ga$ in the Ignatiev model $\cI$ such that $\cI,\ga\models\la n\ra\top,\ [0]\phi$. Let $\gb$ be the uppermost node in which the formula $\la n\ra\top$ is true, that is, the sequence of ordinals $(\gw_n,\gw_{n-1},\dots,\gw,1,0,\dots)$ where $\gw_0=1$ and $\gw_{i+1}=\gw^{\gw_i}$ for $i\geq 0$. We have $\ga R_0\gb$ or $\ga=\gb$. By transitivity, $\cI,\gb\models [0]\phi$. Moreover, for each variable-free $\psi\in\cL(n+1,0)$ and each $\gy$ such that $\gb R_0\gy$, we have $\cI,\gy\models \psi$ iff $\cI,\gy\models \ol{\psi}$. Therefore, $\cI,\gb\models [0]\ol{\phi}$ which implies the claim. \ep

\paragraph{Proof of Theorem~\ref{unif}.} Suppose $\gs$ is a variable-free substitution such that $\GLP_n\vdash \gs(\phi).$ Then $\GLP_n\vdash [0]\gs(\phi)$ and $\GLP\vdash \la n\ra\top \to [0]\gs(\phi).$ If $\GLPz\vdash \la n+1\ra\top \to \neg[0]\phi$ then $\GLPz\vdash \la n\ra\top \to \neg[0]\gs(\phi)$. Therefore, we obtain $\GLPz\vdash [n]\bot$, which is not the case.

In the other direction, suppose $\phi\in\cL(n,m)$. If $\GLPz\nvdash \la n\ra\top \to \neg[0]\phi$ then by Aguilera and Pakhomov there is a variable-free substitution $\gs: \cL(n+1,m)\to \cL(n+1,0)$ such that $\GLP\nvdash \la n\ra\top \to \neg [0]\gs(\phi)$. By Lemma~\ref{bn} this entails $\GLP\nvdash \la n\ra\top \to \neg [0]\ol{\gs(\phi)}$. However, since $[n]$ does not occur in $\phi$, the formula $\ol{\gs(\phi)}$ has the form $\gs'(\phi)$ for some substitution $\gs':\cL(n,m) \to \cL(n,0)$. It follows that $\GLP\nvdash \la n\ra\top \to \neg [0]\gs'(\phi)$ and by Lemma~\ref{vf} $\GLP_n\vdash \gs'(\phi)$. This completes the proof of Theorem~\ref{unif}.

\section{Arithmetical unification and admissibility}
In this section we make a few observations on the arithmetical analogues of the results of this paper. They are based on the standard interpretation of polymodal provability logic, according to which propositional formulas are interpreted as arithmetical sentences, and modalities $[k]$ are interpreted as $k$-provability assertions over $\PA$.

We call an arithmetical formula $\phi$ \emph{$k$-provable}, if $\phi$ is provable from the axioms of Peano arithmetic $\PA$ and all $\Pi_k$-sentences true in the standard model of arithmetic. Let $\Pr_k(x)$ denote a natural formalization of the statement that the formula with the G\"odel number $x$ is $k$-provable (see~\cite{Ign93,Bek05en} for the details).
\renewcommand{\St}{\mathrm{St}}
Let $\St_\PA$ denote the set of sentences in the language of $\PA$.

\bd An \emph{arithmetical realization} is a map $v:\Fm(n,m)\to \St_{\PA}$ such that $v$ commutes with boolean connectives and, for all $\phi$ and $k<n$,
$$v([k]\phi)=\Pr_k(\gn{v(\phi)}).$$
A formula $\phi\in\Fm(n,m)$ is called \emph{arithmetically unifiable} if $\PA\vdash v(\phi)$, for some arithmetical realization $v$. Such a $v$ is then called an \emph{arithmetical unifier} of $\phi$.

An inference rule in the language of $\GLP$ is called \emph{arithmetically admissible}, if every arithmetical unifier of the premisses of the rule unifies its conclusion.
\ed

The analogy with the propositional unification is clear. Notice that both types of unification amount to the same notion applied to two different modal algebras: free GLP-algebras in case of propositional unification, and the provability GLP-algebra of $\PA$ in case of arithmetical unification. Our goal here is to characterize arithmetically unifiable formulas and admissible rules. The key result here is the second arithmetical completeness theorem for $\GLP$ due to G. Japaridze and K. Ignatiev.

\renewcommand{\GLPS}{\mathrm{GLPS}}
Let $\GLPS_n$ denote the closure under the modus ponens and substitution rules of the set of theorems of $\GLP_n$ and all sentences of the form $[k]q\to q$, for $k<n$, that is, $\GLPS_n:=\GLP_n\{[k]q\to q:k<n\}$.

The following theorem is due to G.~Japaridze~\cite{Dzh86} (for a somewhat different arithmetical interpretation) and K.~Ignatiev~\cite{Ign93} (for a more general class of interpretations including this one). See also \cite{Bek11} for a simplified proof.

\bt \label{japcom} For any $n\leq\gw$ and any modal formula $\phi\in\Fm(n,m)$,
\benr \item $\GLP_n\vdash\phi$ iff $\PA\vdash v(\phi)$, for all arithmetical realizations $v$;
\item $\GLPS_n\vdash\phi$ iff $\nat\models v(\phi)$, for all arithmetical realizations $v$.
\eenr
\et

Let $H(\phi)$ denote the formula
$\bigwedge_{i<s}([n_i]\phi_i\to\phi_i)$, where the formulas
$[n_i]\phi_i$, for $i<s$, enumerate all the subformulas of $\phi$ of
the form $[k]\psi$. Japaridze showed that $\GLPS\vdash \phi$ iff $\GLP\vdash H(\phi)\to \phi$, therefore the derivability problem for $\GLPS$ is decidable (see~\cite[Theorem 4]{Bek11}).

As a corollary of Theorem~\ref{japcom} (ii) we obtain the following characterization.

\bt The following statements are equivalent, for any $n<\gw$ and $\phi\in\Fm(n,m)$,\benr
\item $\phi$ is arithmetically unifiable;
\item $\GLPS_n\nvdash \neg [0]\phi$;
\item $\GLP\nvdash \la n\ra\top \to \neg [0]\phi.$
\eenr
\et

\bp\ It is more convenient to prove the equivalence of the negations of the three statements.
Assume (i) does not hold, then, for all arithmetical realizations $v$, $\PA\nvdash v(\phi)$, hence $\nat\models v(\neg [0]\phi)$. Therefore, by Theorem~\ref{japcom} (ii), $\GLPS_n\vdash \neg [0]\phi$. Hence, (ii) does not hold.

Assume (ii) does not hold, that is, $\GLPS_n\vdash \neg [0]\phi$. We prove that $\GLP\vdash
\la n\ra\top \to \neg [0]\phi.$

First, we remark that $\GLP\vdash [n]([k]\psi\to\psi)$, for all $k<n$, by the same argument as in Lemma~\ref{one}. Since $\phi\in \Fm(n,m)$, it follows that $\GLP\vdash [n]H(\phi)$. Therefore, if $\GLPS_n\vdash \neg[0]\phi$, then $\GLP_n \vdash H(\phi)\to \neg[0]\phi$ and $\GLP\vdash [n]\neg[0]\phi$. Then from $\GLP\vdash [0]\phi \to [n][0]\phi$ we obtain $\GLP\vdash [0]\phi\to [n]\bot$, which is a contraposition of the claim.

Assume (iii) does not hold. Let $v$ be any arithmetical realization. Then, since $\GLP$ is arithmetically sound and $v(\la n\ra\top)$ is true, $\nat\models v(\neg[0]\phi)$. Hence, $\PA\nvdash v(\phi)$, that is, $\phi$ is not arithmetically unifiable. \ep

\bt \label{rules}
A rule $\phi_1,\dots,\phi_n/\phi$ in the language of $\GLP$ is arithmetically admissible iff $\GLPS\vdash [0]\phi_1\land\dots\land [0]\phi_n\to [0]\phi$.
\et

\bp\ (if) Assume $v$ is an arithmetical unifier of $\phi_1,\dots,\phi_n$, that is, $\PA\vdash v(\phi_i)$, for all $1\leq i\leq n$. Then $$\nat\models v([0]\phi_1\land\dots\land [0]\phi_n).$$ By the soundness part of Theorem~\ref{japcom} (ii) we conclude that $\nat\models v([0]\phi)$ and $\PA\vdash v(\phi)$, as required.

(only if) Assume $\GLPS\nvdash [0]\phi_1\land\dots\land [0]\phi_n\to [0]\phi$. Then by the completeness part of Theorem~\ref{japcom} (ii) there is an arithmetical realization $v$ such that $\nat\models v([0]\phi_1\land\dots\land [0]\phi_n)$ and $\nat\nmodels v([0]\phi)$. This means that $\PA\vdash v(\phi_i)$, for all $1\leq i\leq n$, and $\PA\nvdash v(\phi)$, that is, $v$ unifies the premisses of the rule, but not the conclusion. \ep

\bcor Arithmetical unification and admissibility problems for $\GLP$ are decidable. \ecor

\medskip
Finally, we make some remarks on the arithmetical analogues of the results in Section~\ref{boxone}. An \emph{arithmetical unifier} of a formula $[1]p$ can be identified with a sentence $A$ such that $\PA\vdash \Pr_1(\gn{A})$. Such formulas were studied in \cite{KolBek19} and called \emph{provably 1-provable}. One of the main results in \cite{KolBek19} is a characterization of provably 1-provable formulas as theorems of the extension of $\PA$ by $\ge_0$-times iterated local reflection schema.

In this result, the ordinal $\ge_0$ reflects the fact that we deal with arithmetical interpretation over $\PA$. A weaker theory $\EA$ as a base theory would yield essentially the same result as in Lemma~\ref{two}: Arithmetical unifiers of $[1]p$ in $\EA$ are theorems of the extension of $\EA$ by the local reflection schema $\Pr_{\EA}(\gn{A})\to A$, for all sentences $A$.

\section{Open questions}
A number of questions related to unification and admissibility in $\GLP$ remain open. We gather some of the most natural ones for the record here.

We note that Ghilardi's notion of \emph{projective formula} makes sense in the context of $\GLP$. Similarly, de Jongh and Visser's notion of \emph{exact formula} equally makes sense. Exact formulas can be defined as the relators of finitely generated subalgebras of a free $\GLP$-algebra. The uniform Craig's interpolation property for $\GLP$~\cite{Sham11} implies that every such finitely generated subalgebra is finitely presented, that is, exact formulas do exist~(see a similar proof for $\GL$ in \cite{AB04}). Projective formulas are the relators of retracts of finitely generated free $\GLP$-algebras, hence every projective formula is exact.

\paragraph{Problem 1.} Are projective formulas the same as exact formulas in $\GLP$?  For $\GL$ both notions are known to be equivalent (though we do not know a published reference).

\paragraph{Problem 2.} Is there a characterization of projective formulas in $\GLP$ in terms of operations on Kripke models or topological models?

For $\GL$ it is well-known that a formula is projective iff it is preserved under the extension of any Kripke model downwards by a new root point (with a suitable evaluation at this point). N. Lukashov established a similar characterization for the polymodal logic J.

\paragraph{Problem 3.} Is there a reasonable description of the set of all unifiers of a given formula in $\GLP$?

Obviously, we cannot give a description in terms of finitely many maximal unifiers. However, we have a hope for some more general kind of characterization.

\paragraph{Problem 4.} Is the admissibility problem for $\GLP$ decidable?

\paragraph{Problem 5.} Is there a reasonable proof system axiomatizing admissible (multi-conclusion) rules in $\GLP$?

\bibliographystyle{plain}
\bibliography{ref-all2.bib}

\end{document}

%% file: macros.tex
\newtheorem{theorem}{Theorem}
\newtheorem{lemma}{Lemma}

\newtheorem{corollary}[lemma]{Corollary}
\newtheorem{proposition}[lemma]{Proposition}

\newtheorem{definition}{Definition}
\newtheorem{example}[lemma]{Example}
\newtheorem{remark}[lemma]{Remark}

\newcommand{\bl}{\begin{lemma}}
\newcommand{\el}{\end{lemma}}
\newcommand{\bt}{\begin{theorem}}
\newcommand{\et}{\end{theorem}}
\newcommand{\bcor}{\begin{corollary}}
\newcommand{\ecor}{\end{corollary}}
\newcommand{\bp}{\begin{proof}\textbf{.}\ }
\newcommand{\ep}{$\boxtimes$ \par\addvspace{\topsep} \end{proof}}
\newcommand{\bpr}{\begin{proposition}}
\newcommand{\epr}{\end{proposition}}
\newcommand{\brem}{\begin{remark} \em}
\newcommand{\erem}{\end{remark}}
\newcommand{\bd}{\begin{definition} \em}
\newcommand{\ed}{\end{definition}}
\newcommand{\bex}{\begin{example} \em
}
\newcommand{\eex}{\end{example}}
\newcommand{\beq}{\begin{equation} }
\newcommand{\eeq}{\end{equation}}

\newcommand{\bi}{\begin{itemize}
  }
\newcommand{\ei}{\end{itemize}}
\newcommand{\ben}{\begin{enumerate} }
\newcommand{\een}{\end{enumerate} }

\newcommand{\refeq}[1]{(\ref{#1})}
\newenvironment{enumr}{

\begin{enumerate}     }{\end{enumerate}

}

\newcommand{\benr}{\begin{enumr}
  }
\newcommand{\eenr}{
\end{enumr}}

\newcommand{\ignore}[1]{}

\newcommand{\al}[1]{\forall #1\:}
\newcommand{\ex}[1]{\exists #1\:}

\newlength{\hilflh}

\newcommand{\naturals}{\mathbb{N}}

\renewcommand{\emptyset}{\varnothing}

\newcommand{\cL}{{\mathcal L}}

\newcommand{\cK}{{\mathcal K}}
\newcommand{\cA}{{\mathcal A}}

\newcommand{\cW}{{\mathcal W}}

\newcommand{\cI}{\mathcal{I}}

\newcommand{\ga}{\alpha}
\newcommand{\gb}{\beta}

\renewcommand{\ge}{\varepsilon}

\newcommand{\gs}{\sigma}
\newcommand{\gy}{\gamma}
\newcommand{\gw}{\omega}

\renewcommand{\phi}{\varphi}

\newcommand{\imp}{\rightarrow}
\newcommand{\eqv}{\leftrightarrow}

\renewcommand{\Pr}{\mathrm{Pr}}

\newcommand{\ol}{\overline}

\newcommand{\GL}{\mathbf{GL}}
\newcommand{\GLP}{\mathbf{GLP}}

\newcommand{\St}{\mbox{\sf St}}

\newcommand{\PA}{\mathsf{PA}}
\newcommand{\EA}{\mathrm{EA}}

\newcommand{\gn}[1]{\ulcorner #1 \urcorner}

\newcommand{\la}{\langle}
\newcommand{\ra}{\rangle}

\renewcommand{\models}{\vDash}      

\newcommand{\Imp}{\Rightarrow}

\newcommand{\nat}{\naturals}

\newcommand{\J}{\mathbf{J}}

\renewcommand{\leq}{\leqslant}
\renewcommand{\geq}{\geqslant}
\newcommand{\GLPS}{\mathbf{GLPS}}


%% file: unification.bbl
\begin{thebibliography}{10}

\bibitem{AgPa24}
J.~Aguilera and F.~Pakhomov.
\newblock The logic of correct models.
\newblock Technical report, 2024.
\newblock arXiv:2402.15382 [math.LO].

\bibitem{AB04}
S.N. Artemov and L.D. Beklemishev.
\newblock Provability logic.
\newblock In D.~Gabbay and F.~Guenthner, editors, {\em {Handbook of
  Philosophical Logic}, 2nd ed.}, volume~13, pages 229--403. Kluwer, Dordrecht,
  2004.

\bibitem{BekJoo12}
L.~Beklemishev and J.~Joosten.
\newblock Problems collected at the {Wormshop 2012 in Barcelona}.
\newblock https://homepage.mi-ras.ru/$\sim$bekl/ Problems/worm$\underline{\
  }$problems.pdf, 2012.

\bibitem{Bek04}
L.D. Beklemishev.
\newblock Provability algebras and proof-theoretic ordinals, {I}.
\newblock {\em {Annals of Pure and Applied Logic}}, 128:103--123, 2004.

\bibitem{Bek05en}
L.D. Beklemishev.
\newblock Reflection principles and provability algebras in formal arithmetic.
\newblock {\em {Russian Mathematical Surveys}}, 60(2):197--268, 2005.
\newblock Russian original: \emph{Uspekhi Matematicheskikh Nauk}, 60(2): 3--78,
  2005.

\bibitem{Bek05a}
L.D. Beklemishev.
\newblock Veblen hierarchy in the context of provability algebras.
\newblock In P.~H\'ajek, L.~Vald\'es-Villanueva, and D.~Westerst{\aa}hl,
  editors, {\em {Logic, Methodology and Philosophy of Science, Proceedings of
  the Twelfth International Congress}}, pages 65--78. Kings College
  Publications, London, 2005.
\newblock Preprint: Logic Group Preprint Series 232, Utrecht University, June
  2004.

\bibitem{Bek10}
L.D. Beklemishev.
\newblock Kripke semantics for provability logic {GLP}.
\newblock {\em {Annals of Pure and Applied Logic}}, 161:756--774, 2010.

\bibitem{Bek11}
L.D. Beklemishev.
\newblock A simplified proof of the arithmetical completeness theorem for the
  provability logic {GLP}.
\newblock {\em Trudy Matematicheskogo Instituta imeni V.A.Steklova},
  274(3):32--40, 2011.
\newblock English translation: \emph{Proceedings of the Steklov Institute of
  Mathematics}, 274(3):25--33, 2011.

\bibitem{Bek17}
L.D. Beklemishev.
\newblock On the reduction property for {GLP-algebras}.
\newblock {\em Doklady: Mathematics}, 95(1):50--54, 2017.

\bibitem{BJV}
L.D. Beklemishev, J.~Joosten, and M.~Vervoort.
\newblock A finitary treatment of the closed fragment of {J}aparidze's
  provability logic.
\newblock {\em Journal of Logic and Computation}, 15(4):447--463, 2005.

\bibitem{Bezh16}
N.~Bezhanishvili, D.~Gabelaia, S.~Ghilardi, and M.~Jibladze.
\newblock Admissible bases via stable canonical rules.
\newblock {\em Studia Logica}, 104(2):317--341, 2016.

\bibitem{Fef60}
S.~Feferman.
\newblock Arithmetization of metamathematics in a general setting.
\newblock {\em Fundamenta Mathematicae}, 49:35--92, 1960.

\bibitem{Ghi99}
S.~Ghilardi.
\newblock Unification in intuitionistic logic.
\newblock {\em The Journal of Symbolic Logic}, 64:859--880, 1999.

\bibitem{Ghi00}
S.~Ghilardi.
\newblock Best solving modal equations.
\newblock {\em Annals of Pure and Applied Logic}, 102:183--198, 2000.

\bibitem{Ica09}
T.F. Icard.
\newblock A topological study of the closed fragment of {GLP}.
\newblock {\em Journal of Logic and Computation}, 21(4):683--696, 2011.

\bibitem{Iem01}
R.~Iemhoff.
\newblock On the admissible rules of intuitionistic propositional logic.
\newblock {\em The Journal of Symbolic Logic}, 66(1):281--294, 2001.

\bibitem{IemMet09}
Rosalie Iemhoff and George Metcalfe.
\newblock Proof theory for admissible rules.
\newblock {\em Annals of Pure and Applied Logic}, 159(1-2):171--186, 2009.

\bibitem{Ign93}
K.N. Ignatiev.
\newblock On strong provability predicates and the associated modal logics.
\newblock {\em The Journal of Symbolic Logic}, 58:249--290, 1993.

\bibitem{Dzh86}
G.K. Japaridze.
\newblock The modal logical means of investigation of provability.
\newblock Ph.D.~{T}hesis in {P}hilosophy, in {R}ussian, Moscow, 1986.

\bibitem{Jera05}
E.~Je{\v r}{\'a}bek.
\newblock Admissible rules of modal logics.
\newblock {\em Journal of Logic and Computation}, 15(4):411--431, 2005.

\bibitem{Jera15}
E.~Je{\v r}{\'a}bek.
\newblock Blending margins: The modal logic {$\mathbf K$} has nullary
  unification type.
\newblock {\em Journal of Logic and Computation}, 25(5):1231--1240, 2015.

\bibitem{KolBek19}
E.~Kolmakov and L.~Beklemishev.
\newblock Axiomatization of provable $n$-provability.
\newblock {\em The Journal of Symbolic Logic}, 84(2):849–869, 2019.

\bibitem{LukBT}
N.~Lukashov.
\newblock Unification problem for the provability logic {GLP}.
\newblock Bachelor's Thesis, National Research University Higher School of
  Economics, Faculty of Mathematics, in Russian, 2023.

\bibitem{Pakh16}
F.N. Pakhomov.
\newblock Linear $\mathrm{GLP}$-algebras and their elementary theories.
\newblock {\em Izvestiya: Mathematics}, 80(6):1159--1199, 2016.

\bibitem{Ryb84}
V.V. Rybakov.
\newblock A criterion for admissibility of rules in the modal system {${\bf
  S4}$} and intuitionistic logic.
\newblock {\em Algebra and Logic}, 23:369--384, 1984.

\bibitem{Ryb89}
V.V. Rybakov.
\newblock On admissibility of the inference rules in the modal system {$G$}.
\newblock In Yu.L. Ershov, editor, {\em Trudy instituta matematiki}, volume~12,
  pages 120--138. Nauka, Novosibirsk, 1989.
\newblock In {R}ussian.

\bibitem{Ryb97}
V.V. Rybakov.
\newblock {\em Admissibility of logical inference rules}.
\newblock Elsevier, Amsterdam, 1997.

\bibitem{Sham11}
D.S. Shamkanov.
\newblock Interpolation properties of provability logics {GL} and {GLP}.
\newblock {\em Trudy Matematicheskogo Instituta imeni V.A. Steklova},
  274(3):329--342, 2011.
\newblock English translation: \emph{Proceedings of the Steklov Institute of
  Mathematics}, 274(3):303--316, 2011.

\bibitem{WoZa08}
F.~Wolter and M.~Zakharyaschev.
\newblock Undecidability of the unification and admissibility problems for
  modal and description logics.
\newblock {\em ACM Transactions on Computational Logic}, 9(4):1--20, 2008.

\end{thebibliography}
